\documentclass[11pt]{article}

\usepackage{amssymb}
\usepackage{mathrsfs}
\usepackage{latexsym}
\usepackage{amsmath}
\usepackage{amsthm}
\usepackage{amsfonts}
\usepackage{url}
\usepackage{amssymb}
\usepackage{bbm}

\usepackage[round]{natbib}

\renewcommand{\vec}[1]{\bs{\mathrm{#1}}}
\newcommand{\hide}[1]{}

\newcommand{\mexp}{\vec{E}}

\newcommand{\E}{\mexp}

\newcommand{\inv}{^{-1}} %

\newcommand{\basicspace}{\calX}
\newcommand{\X}{\basicspace}

\newcommand{\calX}{\mathcal{X}}

\newcommand{\R}{\mathbb{R}}
\newcommand{\N}{\mathbb{N}}

\newcommand{\beq}{\begin{eqnarray*}}
\newcommand{\eeq}{\end{eqnarray*}}
\newcommand{\beqn}{\begin{eqnarray}}
\newcommand{\eeqn}{\end{eqnarray}}

\newcommand{\paren}[1]{\left( #1 \right)}
\newcommand{\sqprn}[1]{\left[ #1 \right]}
\newcommand{\tlprn}[1]{\left\{ #1 \right\}}
\newcommand{\set}[1]{\tlprn{#1}}
\newcommand{\abs}[1]{\left| #1 \right|}
\newcommand{\ceil}[1]{\ensuremath{\left\lceil#1\right\rceil}}
\newcommand{\floor}[1]{\ensuremath{\left\lfloor#1\right\rfloor}}

\newcommand{\bs}{\boldsymbol}

\newcommand{\oo}[1]{\frac{1}{#1}}

\newcommand{\ben}{\begin{enumerate}}
\newcommand{\een}{\end{enumerate}}
\newcommand{\bit}{\begin{itemize}}
\newcommand{\eit}{\end{itemize}}

\newcommand{\del}{\partial}

\newcommand{\dd}[2]{\frac{d#1}{d#2}}
\newcommand{\ddel}[2]{\frac{\del#1}{\del#2}}

\newcommand{\ddelll}[3]{\frac{\del^2#1}{{\del#2}{\del#3}}}

\def\eps{\varepsilon}

\newcommand{\pbig}{\bar\pi}
\newcommand{\psml}{\pi}

 \theoremstyle{plain}
 \newtheorem{thm}{Theorem}
 \theoremstyle{plain}    
 \newtheorem{lem}[thm]{Lemma}
 \theoremstyle{plain}    
 \newtheorem{rem}[thm]{Remark}
 \theoremstyle{plain}    
 \newtheorem{prop}[thm]{Proposition}

\newcommand{\bepf}{\begin{proof}}
\newcommand{\enpf}{\end{proof}}

\title{The Missing Mass Problem}
\author{Daniel Berend and Aryeh Kontorovich}

\begin{document}
\maketitle
\begin{abstract}
We give tight lower and upper bounds on the expected missing mass for distributions over finite and countably infinite spaces. An essential characterization of the extremal distributions is given. We also provide an extension to totally bounded metric spaces that may be of independent interest.
\end{abstract}

\section{Introduction}
\subsection{Background}
\label{sec:bkgnd}

Let $S$ be a countable set endowed with a probability measure $P$. Suppose that $X_1,\ldots,X_t$
are drawn independently from $S$ according to $P$. Define the {\em missing mass} $U_t$ as the following
random variable:
\beqn
\label{eq:Utdef}
U_t &=& P(S\setminus\set{X_1,\ldots,X_t}).
\eeqn
In words, $U_t$ is the total probability mass 
of the elements of $S$ that were not observed in the $t$ samples.

The missing mass
is a quantity of interest 
in almost any application involving sampling from a large discrete set, whether it be 
fish in a pond
or
words in a language corpus.
(Famously, Alan Turing developed what became known as Good-Turing frequency estimators 
\citep{MR0061330}
as part of his work on cracking the Enigma cypher during WWII; see the account in
\citep{MR1807533}).
We note right away that
\beq
\E[U_t] &=& \sum_{s\in S} p_s(1-p_s)^t
\eeq
where $p_s=P(\set{s})$
and that $\E[U_t]\to0$ as $t\to\infty$ 
(the latter follows from Lebesgue's Dominated Convergence Theorem).
Observe also that $U_t\to0$ almost surely as $t\to\infty$;
one way of seeing this is to apply the deviation inequality
\beqn
\label{eq:uconc}
P[
|U_t - \E[U_t]|
\ge 
\eps]&\le& 2e^{-t\eps^2}
\eeqn
of \citet[Theorem 16]{DBLP:journals/jmlr/McAllesterO03}
together with the Borel-Cantelli Lemma.

The topic of interest of this paper is the {rate of decay} of
$\E[U_t]$. For example, when $S$ is finite, we have the trivial estimate
\beq
\E[U_t]&\le& (1-p_{\min})^t
\eeq
where $p_{\min} = \min_{s\in S} p_s$ and we assume without loss of generality
that $P$ has full support. Of course, this bound is not distribution-free, since it depends on $p_{\min}$. Is a distribution-free estimate possible, at least for finite $S$? What about countable $S$? How about lower bounds on the decay rate of the missing mass? These and related questions are investigated in this paper.

\subsection{Related work}
The missing mass problem is 
an unavoidable feature of
density estimation, 
where non-zero density must be assigned to unobserved regions. The process of transferring some of the probability mass from observed points to unobserved ones is known as {\em smoothing}, and 
Laplace's ``add one'' estimator \citep{laplace1814}
appears to be the earliest 
smoother.
Smoothing is now an indispensable component of density estimation; dozens of methods have been proposed for discrete densities alone (cf. 
\citet{MR633417,DBLP:conf/focs/OrlitskySZ03}).

Since a review of smoothing methods is beyond the scope of this paper, we briefly mention the celebrated Good-Turing estimator for the missing mass. 
Given the sample $X=\set{X_1,\ldots,X_t}$, define $Y^{(1)}\subseteq X$ to consist of those $X_i$ that occur exactly once.
The 
Good-Turing missing mass estimator
is given by
\beq
\hat U_t = \oo t \abs{Y^{(1)}};
\eeq
this is the 
proportion
of frequency-one elements.
An attractive feature of this estimator is its diminishing bias:
\beqn
\label{eq:Ubias}
\E[\hat U_t] 
-
\E[U_t] 
&=& \oo t \E[ U_t^{(1)}]
\eeqn
where $U_t^{(1)}=P(Y_t^{(1)})$ is the random variable corresponding to the total mass of the frequency-one items; this variant of Good's theorem is proved in
\citet[Theorem 1]{DBLP:conf/colt/McAllesterS00}.
Additionally, both the missing mass $U_t$ and its
estimate $\hat U_t$ are 
tightly concentrated about their expectations (inequality 
(\ref{eq:uconc}) establishes this for $U_t$; see 
\citet{DBLP:conf/colt/McAllesterS00}
for other deviation estimates).

\section{Main results}
Although (\ref{eq:uconc}) and (\ref{eq:Ubias}) provide a computationally 
efficient estimator of $\E U_t$,
they do not yield any a priori information about the magnitude of the latter.

To state results, it will be convenient to define the {\em plateau length}
$\ell$ 
of a probability distribution $P$ over $\N$:
\beqn
\label{eq:plateau}
\ell(P) = 
\sup_{0<\alpha<1}\abs{\set{i\in \N: \alpha/2\le p_i < \alpha}}
\eeqn
where $p_i=P(\set{i})$. (Note that $\ell(P)=\infty$ is possible.)

Our first 
two
results 
deal with upper and lower bounds on $\E U_t$. We 
use the notation $[n]=\set{1,2,\ldots,n}$
throughout the paper.
\begin{thm}
\label{thm:UB}
The expected missing mass is bounded above as follows:
\bit
\item[(i)]
For 
$n\in\N$ and $S=[n]$,
\beq \E U_t \le 
\left\{
\begin{array}{ll}
e^{-t/n}, \qquad &t\le n, \\
\frac{n}{et}, \qquad &t>n.
\end{array}\right.
\eeq
\item[(ii)]
For 
$S=\N$, 
\beq \E U_t \le  \frac{\ell(P)}{ct},\eeq
where $c$ is a universal constant.
\eit
\end{thm}
\begin{rem}
It is possible to somewhat (but not by much, see Proposition \ref{prop:tight}) 
improve the bound in (i) in some regimes of $n$ and $t$; this will become apparent from our proofs.
The bound $\frac{n}{et}$ holds everywhere, but
is vacuous when $et<n$.
A slightly better bound of $n(1-1/n)^n/t$ was obtained by R. Boppana in a very elegant way, in response to our question
\citep{Boppana}. 
We took an entirely different route, which also 
basically
characterizes the extremal distribution.
\end{rem}

\begin{prop}
\label{prop:tight}
The estimates in Theorem \ref{thm:UB} are essentially tight:
\bit
\item[(i)] For each 
$n\in\N$
and $t>n$, 
there is a distribution 
on $[n]$
such that
\beq
\E U_t \ge c\frac{n-1}{t},
\eeq
where $c$ is an absolute constant.
\item[(ii)] For 
each integer $a>1$,
there is a distribution $P$ over $S=\N$
such that $\ell(P)=a$ and
\beq \E U_t \ge c\frac{a}{t},\qquad t>a,
\eeq
where $c$ is an absolute constant.
\eit
\end{prop}

Furthermore, if we allow distributions with infinite plateau length, then no 
nontrivial
uniform 
(or even pointwise)
bound on $\E U_t$ is possible:
\begin{prop}
\label{prop:rate-lb}
For any sequence $1>r_1>r_2>\ldots$ decreasing to $0$, there is a distribution
on $S=\N$ such that $\E U_t >r_t$ for all $t\ge1$.
\end{prop}

Next, we turn our attention to extremizing distributions for finite $S$. These turn out to exhibit a fairly regular
behavior, with essentially a single phase transition. Since 
$\E U_t$ is a symmetric function of the $\set{p_i}$,
we henceforth assume that
$p_1\le p_2\le\ldots\le p_n$.
In the sequel, the vector $(p_1,\ldots,p_n)$ will be denoted by $\vec p$.

\begin{thm}
\label{thm:extreme}
Let $|S|=n<\infty$. Then
\bit
\item[(i)] 
Every local maximum 
$\vec{p}^*$ 
of
$\E U_t$ 
is 
of the form
\beq p^*_1 = p^*_2 = p^*_3 =\ldots = p^*_{n-1}\le p^*_n \eeq
(that is, $\vec p^*$ consists of one ``heavy'' element and $n-1$ identical ``light'' ones, where
the possibility of
heavy=light is not excluded).
\item[(ii)] 
There exists a threshold 
$\tau=\tau(n)> n$
such that:
\bit
\item[(a)]
For 
$t< \tau$,
there is a unique global maximum
\beq p^*_1 = p^*_2 = p^*_3 =\ldots = p^*_{n-1}=p^*_n=\oo n. \eeq
\item[(b)]
For $t> \tau$,
there is a unique global maximum and it has the form
\beq p^*_1 = p^*_2 = p^*_3 =\ldots = p^*_{n-1}<p^*_n.\eeq
\eit
\item[(iii)] 
As $n\to\infty$,
\beq\tau = n+\sqrt{2n}(1+o(1)).\eeq
\item[(iv)]
For $t
\ge n+\sqrt{2n}
$,
\beqn
\label{eq:e^n/2}
\oo{t+1}
<
p^*_2
<
\oo{t+1}+e^{-\sqrt{n/2}}.
\eeqn
\eit
\end{thm}
\begin{rem}
We have not excluded the possibility that for 
$t=\tau$, both of the distributions described in (ii) attain the maximum.
(Of course, this seems highly improbable.) 
\end{rem}

\section{Proofs}

\begin{lem}
\label{lem:1et}
Consider the function
$f(x)=x(1-x)^t$ 
on the interval $[0,1]$
for an arbitrary fixed $t>0$.
\bit
\item[(i)]
For $t>0$, 
$f$
increases on $(0,1/(t+1))$, decreases on $(1/(t+1),1)$,
and achieves its maximum
at 
$x=1/(t+1)$,
where it is 
bounded above by $\oo{et}$.
\item[(ii)]
The derivative $f'$ decreases on $(0,2/(t+1))$ and increases
on $(2/(t+1),1)$.
\eit
\end{lem}
\bepf
\bit
\item[(i)]
A simple calculation 
shows
that
$g'>0$ on $(0,1/(t+1))$, 
$g'<0$ on $(1/(t+1),1)$,
and
$g'=0$ at 
$x^*=t/(t+1)$. Substituting, we obtain
\beq
g(x^*) &=& \frac{1}{t+1}\paren{1-\oo{t+1}}^t 
= \oo{t}\paren{1-\oo{t+1}}^{t+1} 
< \oo{et}.
\eeq
\item[(ii)]
Routine.
\eit
\enpf

\bepf[of Theorem \ref{thm:UB}]
\bit
\item[(i)] 
It follows from (\ref{eq:bivalent}) below that,
for
$t\le n$,
the expected missing mass is maximized
when
$p^*_1=p^*_2=\ldots=p^*_n=1/n$, yielding the bound
\beq \E U_t = (1-1/n)^t \le e^{-t/n}.\eeq
For general $t$,
an application of Lemma \ref{lem:1et} yields
\beq
\E U_t &=& \sum_{i=1}^n p_i(1-p_i)^t
\le \frac{n}{et}.
\eeq
\item[(ii)] 
Let $P$ be a distribution on $S=\N$. Then
\beq
\E U_t &=& \sum_{i=1}^\infty p_i(1-p_i)^t \\
&=& 
\sum_{i:p_i\ge 1/(t+1)} p_i(1-p_i)^t 
+
\sum_{i:p_i< 1/(t+1)} p_i(1-p_i)^t
\equiv E_1+E_2.
\eeq
By Lemma \ref{lem:1et},
\beqn
E_1 &=& \sum_{i:p_i\ge 1/(t+1)} p_i(1-p_i)^t \nonumber\\
&=& \sum_{j=0}^{\floor{\log_2(t+1)}} \sum_{i: 2^j/(t+1)\le p_i < 2^{j+1}/(t+1)} p_i(1-p_i)^t \nonumber\\
&\le& \sum_{j=0}^{\floor{\log_2(t+1)}} 
\ell(P)
\frac{2^j}{t+1}\paren{1-\frac{2^j}{t+1}}^t \nonumber\\
&=&
\frac{\ell(P)}{t+1}
\sum_{j=0}^{\floor{\log_2(t+1)}}
2^j \paren{1-\frac{2^j}{t+1}}^{(t+1)\cdot t/(t+1)} \nonumber\\
&\le&
\frac{\ell(P)}{t+1}
\sum_{j=0}^{\floor{\log_2(t+1)}}
2^j \exp(-2^jt/(t+1))
\nonumber\\
&\le&
c'\frac{\ell(P)}{t+1}
\sum_{j=0}^{\floor{\log_2(t+1)}}
2^j 
e^{-j}
\nonumber\\
&<&
c'\frac{\ell(P)}{t+1}
\sum_{j=0}^{\infty}
(2/e)^j \nonumber\\
&\le& 
\label{eq:E1bd}
c'' \frac{\ell(P)}{t+1},
\eeqn
for appropriate
absolute constants
$c',c''$.

An analogous argument shows that
\beqn
E_2&\le&
\label{eq:E2bd}
c'''\frac{\ell(P)}{t+1}
.
\eeqn
Combining (\ref{eq:E1bd}) and (\ref{eq:E2bd}), we obtain
the claim.
\eit
\enpf

\bepf[of Proposition \ref{prop:tight}]
\bit
\item[(i)]
Define the distribution $\vec{p}$ by
\beq 
x=p_1 = p_2 =\ldots = p_{n-1} \le p_n=1-(n-1)x
,\eeq
where $x=1/(t+1)$.
Then
\beq
\E U_t &>& \frac{n-1}{t+1}\paren{1-\oo{t+1}}^t 
= \frac{n-1}{t}\paren{1-\frac{1}{t+1}}^{t+1} 
\ge
\frac{8(n-1)}{27t}
.
\eeq
\item[(ii)]
For any $a\in\N$, define $
\vec{p}
$ as follows:
\beq p_1=p_2=\ldots=p_a=\oo{2}; p_{a+1}=\ldots=p_{2a}=\oo{4a};
\ldots; p_{ka+1}=\ldots=p_{(k+1)a}=\oo{2^ka};\ldots. \eeq
Then, denoting $\kappa=\ceil{\log_2(t/a)}$, we have for $t>a$:
\beq
\E U_t &=& \sum_{k=1}^\infty \oo{2^k}\paren{1-\oo{2^ka}}^t
> \oo{2^\kappa}\paren{1-\oo{2^\kappa a}}^t
= \oo2\cdot \oo{2^{\kappa-1}}\paren{1-\oo{2^\kappa a}}^t
> \frac{a}{2t} \paren{1-\oo{t}}^t
\\
&\ge& 
\frac{4a}{27 t}
.
\eeq
\eit
\enpf

\bepf[of Proposition \ref{prop:rate-lb}]
Let $1>r_1>r_2>\ldots$ be a sequence decreasing to $0$. Observe that
\beq
\E U_t &=& \sum_{i=1}^n p_i(1-p_i)^t 
\ge \sum_{i:p_i<1/t^2} p_i\paren{1-\oo{t^2}}^t
= \paren{1-\oo{t^2}}^t \sum_{i:p_i<1/t^2} p_i.
\eeq
Select
$\tau>10$ such that
$r_{\tau}<0.9$.
Then
we can choose
$(p_i)$ so that
\beqn
\label{eq:tau}
\paren{1-\oo{t^2}}^t \sum_{i:p_i<1/t^2} p_i
> r_t,
\qquad t\ge\tau.
\eeqn
Indeed, 
$(1-1/t^2)^t>0.9$
for $t\ge\tau>10$. Thus, for $t=\tau$, (\ref{eq:tau}) is satisfied by any 
$(p_i)$ with
$p_i<1/\tau^2$ for all $i\in\N$. 
For each $t>\tau$, choose a finite sequence $(p_{it})$ such that
$p_{it}<1/(t+1)^2$ for each $i$ and
\beq
\sum p_{it}
&=& r_{t+1}-r_t.
\eeq
Let $\vec{p}$ be the distribution obtained by concatenating
all the sequences $(p_{it})_{t>\tau}$ and the number $1-r_\tau$.
To prove the claim for $t<\tau$, let us define the following 
``doubling operator'' on distributions:
\beq D((p_1,p_2,\ldots)) = (p_1/2,p_1/2,p_2/2,p_2/2,\ldots).\eeq
It is straightforward to verify that for all 
distributions $\vec{p}$ and all
$t\in\N$,
\beq 
\E_{D^k\vec{p}} U_t 
\nearrow 1
\text{ as }k\to\infty
\eeq
(where 
the subscript of $\E$
specifies the distribution under which the 
expectation is taken).
Thus, if $\vec{p}$ is a distribution that satisfies (\ref{eq:tau}) for all $t\ge\tau$, there is some
finite $k$ such that $D^k\vec{p}$ makes the proposition hold.
\enpf

\bepf[of Theorem \ref{thm:extreme}]
\bit
\item[(i)]
For $n,t\in\N$, define $F:[0,1]^n\to\R$ by
\beq F(\vec{x})=\sum_{i=1}^n x_i(1-x_i)^t
=
\sum_{i=1}^n f(x_i)
\eeq
(where $f(x)=x(1-x)^t$).
An elementary application of Lagrange multipliers shows that, 
under the 
constraint $\sum_{i=1}^n x_i=1$,
a necessary condition for an extremum is
\beq \ddel{F}{x_i}=\ddel{F}{x_j} ,
\qquad i,j\in[n].
\eeq
Lemma \ref{lem:1et}
leaves two possibilities for an extreme point $\vec{p^*}$:
either
all the $p_i^*$ take the value $1/n$ 
(we call such distributions {\em univalent})
or
the $p_i^*$ take
two values
$\psml<\pbig$, with $f'(\psml)=f'(\pbig)<0$
(we call such distributions {\em bivalent}).

In the bivalent case, we have, without loss of generality,
\beq p_1^*=p_2^*=\ldots=p_k^*=\psml<\pbig=p_{k+1}^*=p_{k+2}^*=\ldots=p_n^*\eeq
for some $1<k<n$.
Define the Lagrangian
\beq L(\vec{x},\lambda) = F(\vec{x})+\lambda (g(\vec{x})-1),\eeq
where $g(\vec{x})=\sum_{i=1}^n x_i$ and the 
associated $(n+1)\times(n+1)$ matrix
$H=H(\vec{x},\lambda)$, where
\beq
H_{ij}&=&
\left\{
\begin{array}{ll}
\ddelll{L}{x_i}{x_j}=t(x_i(t+1)-2)(1-x_i)^{t-2},\qquad & i=j\le n, \\
\ddelll{L}{x_i}{x_j}=0,\qquad & i\neq j\le n, \\
\ddel{g}{x_i}=1,\qquad & i\le n,j=n+1, \\
\ddel{g}{x_j}=1,\qquad & j\le n,i=n+1, \\
0,\qquad \qquad & i=j=n+1.
\end{array}\right.
\eeq

Suppose $k\le n-2$ and consider the 
$3\times3$ lower right
submatrix
\beq
B=B(\vec{x}) &=& \paren{
\begin{array}{ccc}
t(x_{n-1}(t+1)-2)(1-x_{n-1})^{t-2} & 0 & 1\\
0& t(x_{n}(t+1)-2)(1-x_{n})^{t-2} & 1\\
1&1&0
\end{array}
};
\eeq
note that our assumption on $k$ forces $B_{11}=B_{22}$.
The second-order 
necessary
condition for $\vec{p}^*$ to be a local maximum is that
a sequence of
bordered Hessians, including $\det(B(\vec{p}^*))$, be 
nonnegative.
Since 
$f'(\psml)=f'(\pbig)<0$
and $f'$ decreases on $(0,2/(t+1))$ and increases on $(2/(t+1),1)$,
it follows
that 
$\psml<2/(t+1)$ and
$\pbig>2/(t+1)$. 
This implies that
$B_{11}=B_{22}>0$.
Denoting this common value by $b$, we have
\beq \det(B)=-2b<0.\eeq
The necessary condition is thus violated, leaving two possibilities:
the univalent case $p_i^*\equiv 1/n$ 
and the bivalent case with $k=n-1$:
\beqn
\label{eq:bivalent}
 \oo{t+1}<p_1^*=p_2^*=\ldots=p_{n-1}^*<\frac{2}{t+1}<p_n.
\eeqn

\item[(ii)]
We saw above
that $\E U_t$ is always maximized by a 
distribution $\vec p$ 
of the form
\beq 
x=p_1 = p_2 =\ldots = p_{n-1} \le p_n=1-(n-1)x.
\eeq
For distributions of this form, we have
\beqn
\label{eq:gdef}
\E U_t = G_t(x) = (n-1)x(1-x)^t + (1-(n-1)x)((n-1)x)^t,
\eeqn
where $G_t$ is defined on $[0,1/n]$.
Note that $x=1/n$ corresponds to the univalent (uniform) distribution, while
$x<1/n$ corresponds to a bivalent distribution.

We claim the existence of a function $\tau:\N\to\N$ such that\\
(a) for $t<\tau(n)$, $G_t$ 
has the unique maximizer $x^*=1/n$\\
(b) for $t>\tau(n)$, $G_t$ 
has the unique maximizer $x^*<1/n$.

(In principle, it may be possible 
for $G_t$ to have
two distinct maxima on $[0,1/n]$
for $t=\tau(n)$,
but this is rather implausible.)

For $t\le n$, (\ref{eq:bivalent}) implies that 
$G_t$ 
has the unique maximizer
$x^*=1/n$;
this shows that $\tau(n)> n$ (if the
function $\tau$ described in (a) and (b) exists at all). 

Now
define the function $R_t(x)=G_t(x)/G_t(1/n)$.\hide{
We begin with a monotonicity result:
if $\E U_t$ is maximized by a bivalent distribution
for some $t\in\N$, 
then this will also be the case for $t+1$. To 
see this,
define
\beqn
\label{eq:gdef}
G_t(x) = (n-1)x(1-x)^t + (1-(n-1)x)((n-1)x)^t. 
\eeqn
For $t\le n$, (\ref{eq:bivalent}) implies that $G_t$ is maximized 
only
at $x=1/n$;
this shows that $\tau(n)> n$ (if the threshold $\tau(n)$ indeed exists).}
Then
\beqn
\label{eq:Rtdef}
R_t(x)&=&
(n-1)x\paren{\frac{1-x}{1-1/n}}^t 
+ 
(1-(n-1)x)\paren{\frac{(n-1)x}{1-1/n}}^t. 
\eeqn
For $x<1/n$, 
the first term 
on the right-hand side of (\ref{eq:Rtdef})
grows exponentially 
with $t$,
and certainly $R_t(x)>1$ is achieved for some finite $t$.
But 
this
means that $G_t$ 
has a unique maximum
at some $x<1/n$
and so
any function
$\tau(n)$
satisfying (a) and (b) must be
finite
for all $n$.

Suppose that $t$ is such that $G_t$ 
achieves a maximum at $x<1/n$.
It follows from
the uniqueness proof below and from
(\ref{eq:x^*}) 
that the maximizer
$x^*$ of $G_t$ is contained in the interval 
$I_t=(1/(t+1),1/t)$. 
We claim that $R_t(x)\ge1$ implies $R_{t+1}(x)>1$
for all $x\in I_{t}$; from here, 
the existence of $\tau$ satisfying
(a) and (b) follows immediately.
Treating $t$ as a continuous variable,
we have
\beq
\dd{R_t(x)}{t} 
&=&
(n-1)x
\log\paren{\frac{1-x}{1-\oo n}}
\paren{\frac{1-x}{1-\oo n}}^t 
+ 
(1-(n-1)x)
\log\paren{nx}
\paren{nx}^t.
\eeq
We establish
the
monotonicity claim
by showing that
\beqn
\label{eq:R'}
\dd{R_t(x)}{t} 
&>& 0,\qquad t\ge n+1, ~ x\in I_t
.
\eeqn
Indeed, appealing to the inequalities 
$\xi^t\log\xi\ge-\oo{et}$ 
for $0<\xi<1$
and
$\xi^t\log\xi\ge\xi-1$
for $\xi\ge1$
(checked by elementary calculus), 
we see that the inequality
\beq
(n-1)x
\paren{\frac{1-x}{1-1/n}-1}
>
\frac{1-(n-1)x}{et}.
\eeq
is even stronger than
(\ref{eq:R'}).
The latter will hold as long as
\beq
-entx^2+(et+n-1)x-1>0,
\eeq
and it suffices to verify the inequality at 
the endpoints 
$1/(t+1)$ and $1/t$
of
$I_t$,
which is straightforward.
This 
proves the existence of $\tau$ as claimed in (a) and (b).

Uniqueness is established by noting 
(again, via elementary though rather tedious calculus) 
that $G_t'$
vanishes at $x=1/n$ and at not more than two points in the interval
$[1/(t+1),1/n)$, and is strictly positive on $[0,1/(t+1)$.

\item[(iii)] 
For $n\in\N$ and $t\neq\tau(n)$,
let $x^*$ be the maximizer of 
the function $G$ defined in (\ref{eq:gdef}),
where we have dropped the subscript $t$. 
A sufficient condition for
$ G(x^*)>\paren{1-\oo n}^t$
to hold is 
$ G\paren{\oo{t}}>\paren{1-\oo n}^t.$\hide{

}
The latter, in turn, will hold as long as 
$(n-1)(1-1/t)^t/t > (1-1/n)^t$.
We will show that the latter inequality holds for large $n$, if $t=n+\sqrt{2n}$.
To this end,
define the function
\beq
g(n) = \frac{n-1}{n+\sqrt{2n}}\paren{1-\oo{n+\sqrt{2n}}}^{n+\sqrt{2n}}
-
\paren{1-\oo n}^{n+\sqrt{2n}}.
\eeq
For $\nu\in(0,1)$, define $\tilde g(\nu)=g(1/\nu)$ and expand it about $\nu=0$:
\beq
\tilde g(\nu) = \frac{\sqrt 2}{3e}\nu^{3/2} + O(\nu^2).
\eeq
Since for this choice of $t$ we have $G(1/t)-(1-1/n)^t=\Omega_+(n^{-3/2})$, 
it follows
that
\beqn
\label{eq:tauub}
\tau(n) \le n+\sqrt{2n},
\qquad n \gg 1
.
\eeqn
To get a lower bound on $\tau$, we estimate $G(x^*)$ from above:
\beq
G(x^*) &\le& 
\frac{n-1}{t+1}\paren{1-\oo{t+1}}^t
+ \paren{1-\frac{n-1}{t}}\paren{\frac{n-1}{t}}^t
=: \bar G_t(n)
.
\eeq
Now let
\beqn
Q_t(n) 
&=& 
\frac{\bar G_t(n)}{(1-1/n)^t}\nonumber\\
&=& 
\label{eq:Qtn}
\frac{n-1}{t+1}\paren{ \frac{1-1/(t+1)}{1-1/n} }^t
+ \paren{1-\frac{n-1}{t}} \paren{ \frac{(n-1)/t }{1-1/n}   }^t
\eeqn
and note that $Q_t(n)<1$ implies $t<\tau(n)$.
Let us put $t=n+(1-\eps)\sqrt{2n}$ for some $0<\eps<1$, and 
observe
that 
for this choice of $t$,
the second term 
on the right-hand side of (\ref{eq:Qtn})
is negligible:
\beq
\paren{1-\frac{n-1}{t}} \paren{ \frac{(n-1)/t }{1-1/n}   }^t
&<&
\paren{ \frac{(n-1)/t }{1-1/n}   }^t 
= 
\paren{ \frac{n}{t}   }^t 
= 
\paren{ 1-\frac{t-n}{t}   }^t\\
&<&
\exp(-(t-n))
=
\exp({-(1-\eps)\sqrt{2n}}).
\eeq
Let us now
examine the asymptotic behavior of the first term in (\ref{eq:Qtn}).
To this end,
define the functions
\beq
h(n) = 
\frac{n-1}{t+1}\paren{ \frac{1-1/(t+1)}{1-1/n} }^t
\eeq
and $\tilde h(\nu)=h(1/\nu)$, and expand about $\nu=0$:
\beq
\hide{
}
\tilde h(\nu) = 1 - (2-\eps)\eps \nu + O(\nu^{3/2}).
\eeq
Hence, for this choice of $t$, we have
\beq
\frac{G(x^*)}{(1-1/n)^t} \le 1-(2-\eps)\eps n\inv+O(n^{-3/2})
+\exp({-(1-\eps)\sqrt{2n}}) <1
\eeq
for sufficiently large $n$.
This, combined with (\ref{eq:tauub}), implies
\beq
\tau(n)=n+(1+o(1))\sqrt{2n}.
\eeq

\hide{
The latter, in turn, will hold as long as
\beqn
\label{eq:ln<}
\ln(n-1)-\ln(t+1)+t\ln\paren{1-\oo{t+1}} > t\ln\paren{1-\oo n}.
\eeqn
Recalling that 
\beqn
\label{eq:logbds}
-x<\ln(1-x)<-x-x^2/(2-2x), \qquad |x|<1,
\eeqn
we conclude that (\ref{eq:ln<}) is implied by
\beq
\label{eq:tn}
\ln(n-1)-\ln(t+1)
-\frac{2t+1}{2t+2}
> 
-\frac{t}{n}-\frac{t}{2n^2}
.
\eeq
Putting $t=n+h$, we 
apply (\ref{eq:logbds}) again to obtain the stronger inequality
\beq
\label{eq:tnh}
2hn(n-1)(n+h+1)+t(n-1)(n+h+1)>2n^2(n-1)h+3n^2(n-1)+(h+2)^2n^2.
\eeq
Keeping only the highest order terms, we have that
the latter
holds for $h=\sqrt{2n}(1+o(1))$.
}

\item[(iv)]
We claim that
in the bivalent case,
the maximizer $x^*$ of $G$
satisfies
\beqn
\label{eq:x^*}
\oo{t+1} < x^* < \oo{t}.
\eeqn
The first inequality follows from
(\ref{eq:bivalent}). To verify the second inequality it suffices to show that $G'(1/t)<0$.
Indeed,
\beq
\frac{G'(1/t)}{t(n-1)} &<& -\paren{\oo t-\oo{t-1}}\paren{1-\oo t}^{t-1}+\paren{\frac{n-1}{t}}^{t-1} \\
&=& -\oo{t(t+1)}\paren{1-\oo t}^{t-1}+\paren{\frac{n-1}{t}}^{t-1} \\
&<& \paren{\frac{n}{t}}^{t-1}-\oo{t(t+1)(1-1/t)}\paren{1-\oo t}^t \\
&<& \paren{\frac{n}{t}}^{t-1}-\oo{e(t+1)(t-1)} 
\le \paren{\frac{n}{t}}^{t-1}-\oo{et^2}\\
&\le& \frac{n^2}{t^2}\paren{\frac{n}{t}}^{t-3}-\oo{et^2}\\
&\le& \oo{t^2}\sqprn{ \paren{\frac{n}{n+\sqrt{2n}}}^{n-3}-\oo e} <0.
\eeq

To establish 
our claim,
we seek
a small $\delta>0$ such that $G'(1/(t+1)+\delta)<0$;
any such $\delta$ will yield the bound
$ \oo{t+1}<x^*<\oo{t+1}+\delta.$
Putting $p=\oo{t+1}+\delta$, we have
\beq
\frac{G'(p)}{n-1} 
&=&
-\delta(t+1)(1-p)^{t-1}+(t-n+1-(n-1)\delta(t+1))((n-1)p)^{t-1} \\
&<&
-\delta(t+1)(1-p)^{t-1}+(t-n+1)((n-1)p)^{t-1} \\
&<&-\delta t(1-p)^{t-1}+t((n-1)p)^{t-1}.
\eeq
Thus, to ascertain that $G'(p)<0$ it suffices to show that
\beq
\delta &>& \paren{\frac{(n-1)p}{1-p}}^{t-1}.
\eeq
It follows from (\ref{eq:x^*}) that we may take $p<1/t$, and hence
\beq
\paren{\frac{(n-1)p}{1-p}}^{t-1}
<
\paren{\frac{(n-1)/t}{1-1/(t+1)}}^{t-1}
=
\paren{\frac{(n-1)(t+1)}{t^2}}^{t-1}
<
\paren{\frac{n}{t}}^{t-1}.
\eeq
Our assumption that $t\ge n+\sqrt{2n}$ implies
\beq
\paren{\frac{n}{t}}^{t-1} \le
\paren{\frac{n}{n+\sqrt{2n}}}^{t-1} 
\le \exp(-\sqrt{n/2}).
\eeq
Thus, we may take $\delta=e^{-\sqrt{n/2}}$.
\eit
\enpf

\section{Application: missing mass in metric spaces
}
If $P$ is a 
nondegenerate
continuous distribution, then the missing mass as defined in (\ref{eq:Utdef}) is trivially $1$ for all $t\in\N$.
To define a nontrivial extension of this notion to continuous spaces,\footnote{
This 
problem
is of interest in anomaly detection applications
\citep{KHM2011}.
} let us start with a metric probability space $(\X,P,d)$,
whose $\sigma$-field is induced by the metric topology.
For $x\in\X$, let
$B_\eps(x)$ to be the $\eps$-ball about $x$:
$ B_\eps(x) = \set{y\in\X: d(x,y)\le\eps}.$
For $S\subset\X$, define its $\eps$-envelope, $S_\eps$, to be
\beq S_\eps = \bigcup_{x\in S} B_\eps(x) .\eeq
For $\eps>0$, define the $\eps$-covering number, $N(\eps)$, of $\X$
as the minimal cardinality of a set $E\subset\X$ such that
$\X = E_\eps$. A space is 
{\em totally bounded} if $N(\eps)<\infty$ for all $\eps>0$.
Define the
{\em $\eps$-missing mass}
of the sample $S=\set{X_1,\ldots,X_t}$
as the random variable
\beqn
\label{eq:Un}
U_t(\eps) = P(\X\setminus S_\eps).
\eeqn
The expected $\eps$-missing mass of totally bounded spaces is controlled via the covering numbers:
\begin{thm}
In a 
totally bounded
metric probability space $(\X,P,d)$,
\beq \E U_t(\eps) \le \frac{N(\eps)}{et}.\eeq
\end{thm}
\bepf
For a fixed $\eps>0$, 
let $\set{e_1,e_2,\ldots,e_n}$ be an $\eps$-net for $\X$.
For $i=1,\ldots,n$, put $p_i=P(B_\eps(e_i))$; note that $\sum p_i\ge1$.
Then, invoking Lemma \ref{lem:1et}, we have
\beq
\E U_t(\eps) &\le& \sum_{i=1}^n p_i(1-p_i)^t 
\le
\frac{n}{et}.
\eeq
\enpf

\section*{Acknowledgments}
We thank Antonio Cuevas
and Larry Wasserman for helpful correspondence.

\end{document}